\documentclass[11pt]{amsart}
\usepackage[latin1]{inputenc}
\usepackage[T1]{fontenc}
\usepackage[right=2cm,left=2cm,height=23cm]{geometry}

\usepackage{amssymb}
\usepackage{amsmath}
\usepackage{amsthm}

\usepackage{graphicx}
\usepackage{pst-all}
\usepackage{auto-pst-pdf}

\usepackage{xcolor}
\definecolor{bleufonce}{rgb}{0,0,0.2}
\definecolor{couleur1}{rgb}{0.9,0,0}
\definecolor{couleur}{rgb}{0,0,0.5}

\usepackage{hyperref}
\hypersetup{
colorlinks=true,
linkcolor=bleufonce,urlcolor=bleufonce,
citecolor=bleufonce,filecolor=bleufonce}

\newtheorem{FCT}{Fact\!\!}
\newtheorem{PROP}{Proposition\!\!}
\newtheorem{COR}{Corollary\!\!}

\newtheorem{THEA}{Theorem 1\!\!}
\newtheorem{THEB}{Theorem 2\!\!}
\newtheorem{THEC}{Theorem 3\!\!}
\newtheorem{THED}{Theorem 4\!\!}

\newtheoremstyle{Remarques}
  {1pt}
  {1pt}
  {\normalfont}
  {}
  {\bfseries}
  {.}
  {.5em}
  {}
\theoremstyle{Remarques}

\newcommand{\eqnsection}{
\renewcommand{\theequation}{\thesection.\arabic{equation}}
    \makeatletter
    \csname  @addtoreset\endcsname{equation}{section}
    \makeatother}
\eqnsection

\def\a{\alpha}
\def\B{\beta}
\def\Bab{\B_{a,b}}
\def\Babs{\B_{a,b}^{-s}}
 
\def\EE{{\mathbb{E}}} 
\def\ELP{{\mathcal{E}}} 
\def\Ga{\Gamma}
\def\GG{{\mathcal{G}}} 
\def\HH{{\mathcal{H}}} 
\def\II{{\mathcal{I}}} 
\def\MM{{\mathcal{M}}} 

\def\ii{\mathrm{i}}
\def\lbd{\lambda}
\def\lacc{\left\{}
\def\lcr{\left[}
\def\lp{\left(}

\def\racc{\right\}}
\def\rp{\right)}
\def\rcr{\right]}

\def\SS{{\mathcal{S}}} 
\def\U{{\bf U}}
\def\Un{{\bf 1}}

\def\claw{\stackrel{d}{\longrightarrow}}
\def\elaw{\stackrel{d}{=}}
\def\qed{\hfill$\square$}

\newcommand{\C}{\mathbb C}
\newcommand{\R}{\mathbb R}

\newcommand{\E}[1]{\mathbf{E}\left(#1\right)}
\newcommand{\la}{\lambda}
\newcommand{\e}{\varepsilon}
\newcommand{\ex}[1]{\exp\left(#1\right)}

\newcommand{\ga}{\gamma}
\newcommand{\al}{\alpha}
\newcommand{\NN}{\mathbb N}

\renewcommand{\Im}{\mathrm{Im}}
\renewcommand{\Re}{\mathrm{Re}}
               
\title{On the infinite divisibility of inverse Beta distributions}

\author[Pierre Bosch]{Pierre Bosch}

\address{Laboratoire Paul Painlev\'e, Universit\'e Lille 1, Cit\'e Scientifique, F-59655 Villeneuve d'Ascq Cedex. {\em Email} : {\tt pierre.bosch@ed.univ-lille1.fr}}

\author[Thomas Simon]{Thomas Simon}

\address{Laboratoire Paul Painlev\'e, Universit\'e Lille 1, Cit\'e Scientifique, F-59655 Villeneuve d'Ascq Cedex. Laboratoire de physique th\'eorique et mod\`eles statistiques, Universit\'e  Paris Sud, B\^atiment 100, F-91405 Orsay Cedex. {\em Email} : {\tt simon@math.univ-lille1.fr}}



\begin{document}

\begin{abstract} We show that all negative powers $\B_{a,b}^{-s}$ of the Beta distribution are infinitely divisible. The case $b\le 1$ follows by complete monotonicity, the case $b > 1, s\ge 1$ by hyperbolically complete monotonicity and the case $b > 1, s < 1$ by a L\'evy perpetuity argument involving the hypergeometric series. We also observe that $\B_{a,b}^{-s}$ is self-decomposable if and only if $2a + b + s + bs \ge 1,$ and that in this case it is not necessarily a generalized Gamma convolution. On the other hand, we prove that all negative powers of the Gamma distribution are generalized Gamma convolutions, answering to a recent question of L. Bondesson.
\end{abstract}

\maketitle
\section{Introduction and statement of the results}

This paper is a sequel to our previous article \cite{NOUS}, where we have established the infinite divisibility of all negative powers $\ga_a^{-s}$ of the Gamma distribution $\ga_a,$ with density
$$\frac{x^{a-1}}{\Ga(a)}\,e^{-x}\,\Un_{(0,+\infty)}(x).$$
More precisely, in \cite{NOUS} we had completed the already known situation $s\ge 1$ by an argument involving the exponential functional of a spectrally negative L\'evy process, valid in the case $s \le 1$ and in this case only. We consider here the same problem for the Beta distribution $\Bab$, with density
$$\frac{\Ga(a+b)}{\Ga(a)\Ga(b)}\,x^{a -1}(1-x)^{b-1}\,\Un_{(0,1)}(x).$$
Recall that a positive random variable $X$ is infinitely divisible, which we will denote by $X\in\II,$  if and only if its Laplace exponent $\varphi(\lbd) = -\log \EE[e^{-\lbd X}]$ is a Bernstein function, viz.
$$\varphi (\lbd)\; = \; a\lbd\; +\; \int_0^\infty (1- e^{-\lbd x}) \, \nu (dx)$$
with $a \ge 0$ (the drift coefficient) and $\nu$ a non-negative measure on $(0, +\infty)$ (the L\'evy measure) whose integral along $1\wedge x$ is finite. When $\nu$ is absolutely continuous with a density of the type $x^{-1} k(x)$ for some non-increasing function $k$, this means that $X$ is self-decomposable ($X\in\SS$ for short) in other words that for all $c\in (0,1)$ there is a decomposition
$$X\; \elaw \; c X\; +\; X_c,$$
where $X_c$ is an independent random variable. The distribution of $X$ is called a generalized Gamma convolution, which we will denote by $X\in\GG,$ when $X$ is self-decomposable and its above spectral function $k$ is completely monotone (CM). From the probabilistic point of view, the $\GG-$property means that $X$ can be decomposed into the finite or infinite independent sum of Gamma random variables, possibly with different parameters. From the analytical viewpoint, the complete monotonocity of $k$ means that the Laplace exponent of $X$ is a Thorin-Bernstein function. We refer to \cite{Bd,S,SVH} for various accounts on infinite divisibility, self-decomposability and generalized Gamma convolutions, and to \cite{SSV} for a recent monograph devoted to Bernstein functions. 

The positive powers of $\Bab$ cannot be infinitely divisible because of their bounded support - see Theorem 24.3 in \cite{S}. On the other hand, it is well-known that $-\log(\Bab)\in\II$ with an explicit L\'evy measure - see e.g. Example VI.12.21 in \cite{SVH} and the proof of Theorem 1 thereafter for further details. In the present paper, it will be shown that all negative powers of $\Bab$ belong to $\II.$ This allows to retrieve the $\II-$property for 
\begin{equation}
\label{logB}
-\log(\Bab)\; =\; \lim_{s\to 0+}s^{-1}(\Babs -1),
\end{equation}
and also for all the negative powers $\ga_{a}^{-s}$ in view of the convergence in law
\begin{equation}
\label{BG}
 b^{-s}\Babs\; \claw\; \ga_a^{-s} \qquad\mbox{as $b\to +\infty.$}
 \end{equation}
The infinite divisibility of $\Bab^{-s}$ is equivalent to that of $\Bab^{-s}-1,$ whose support is $\R^+.$ In the case $s=1$, the latter random variable is known as the Beta random variable of the second kind, and its infinite divisibility appears in the list of examples of Appendix B in \cite{SVH}. So far, the problem of infinite divisibility for other negative powers of the Beta distribution seems to have escaped investigation. Having independent interest, these random variables appear as multiplicative factors because of their moments of the Gamma type \cite{J}. For instance, they are connected to real stable densities via the Kanter random variable - see (2.4) and (7.1) in \cite{SF}. It was conjectured in \cite{JS} that $\Babs\in\GG$ for all $s \ge 1$ - see Conjecture 2 therein. 

It does not seem possible to express the Laplace exponent of $\Babs$ or $\Babs -1$ in a sufficiently explicit way in order to show that it is a Bernstein function. We will hence proceed via different methods, characterizing properties which are more informative than infinite divisibility or self-decomposability. Let us first introduce the class $\MM$ of positive random variables having a CM density on $(0, +\infty).$ This class is included in $\II$ by Goldie's criterion - see e.g. Theorem 51.6 in \cite{S}.

\begin{THEA} One has $\beta_{a,b}^{-s}-1\in\MM$ if and only if $b\le 1.$
\end{THEA}

Second, let us consider the class $\HH$ of positive random variables having a hyperbolically completely monotone (HCM) density on $(0, +\infty).$ A function $f : (0, +\infty) \to (0, +\infty)$ is said to be HCM if for every $u > 0$ the function $f(uv)f(u/v)$ is CM in the variable $v + 1/v.$ The fact that $\HH\subset\GG$ is a key-result in infinite divisibility on $\R^+$ with explicit density but without explicit Laplace transform, allowing to prove the $\II-$property for many classical positive distributions. We refer to Chapters 4-6 in \cite{Bd} for a complete account, and to Theorem 5.2.1 therein for a proof of the inclusion $\HH\subset\GG.$
 
\begin{THEB} One has $\beta_{a,b}^{-s}-1\in\HH$ if and only if one of the three following conditions is verified.
\begin{enumerate}
\item $b\wedge s>1.$ 
\item $b=1$ or $s = 1.$
\item $b <1,s\in[1/2,1)$ and $a+b+s\ge1.$
\end{enumerate}
\end{THEB}

\bigskip

\begin{minipage}{0.4\linewidth}
\centering
\psset{unit=1.2cm,algebraic=true,linewidth=0.6pt,
,arrowsize=4pt 2,arrowinset=0.2}

\begin{pspicture*}(-1,-1)(4,4)

\rput[t](3.6,-0.28){$s$}
\rput[r](-0.2,3.6){$b$}

\pspolygon[hatchcolor=couleur1,linecolor=couleur1,
fillstyle=hlines,hatchangle=45,hatchsep=0.1](1,1)(6,1)(6,10)(1,10)

\psline[linecolor=couleur1](0,1)(1,1)
\psline[linecolor=couleur1](1,0)(1,1)

\pspolygon[hatchcolor=couleur1,linecolor=couleur1,fillstyle=hlines,hatchangle=45,hatchsep=0.1](0.5,1)(1,1)(1,0)(0.8,0)(0.5,0.3)

\psaxes[Dx=1,Dy=1,ticksize=-2pt 2pt]{->}(0,0)(-0.9,-0.9)(4,4)

\end{pspicture*}\\
The case $a<1/2.$
\end{minipage}
\begin{minipage}{0.05\linewidth}
$ $
\end{minipage}
\begin{minipage}{0.4\linewidth}
\centering
\psset{unit=1.2cm,algebraic=true,linewidth=0.6pt,
,arrowsize=4pt 2,arrowinset=0.2}

\begin{pspicture*}(-1,-1)(4,4)
\rput[t](3.6,-0.28){$s$}
\rput[r](-0.2,3.6){$b$}

\pspolygon[hatchcolor=couleur1,linecolor=couleur1,
fillstyle=hlines,hatchangle=45,hatchsep=0.1](1,1)(6,1)(6,10)(1,10)

\psline[linecolor=couleur1](0,1)(1,1)
\psline[linecolor=couleur1](1,0)(1,1)

\pspolygon[hatchcolor=couleur1,linecolor=couleur1,fillstyle=hlines,hatchangle=45,hatchsep=0.1](0.5,1)(1,1)(1,0)(0.5,0)

\psaxes[Dx=1,Dy=1,ticksize=-2pt 2pt]{->}(0,0)(-0.9,-0.9)(4,4)
\end{pspicture*}\\
The case $a\ge1/2.$
\end{minipage}
\vspace{0.6cm}

So far, we can deduce that $\Babs\in\II$ if $b\le 1$ or $s\ge 1.$ In order to handle the situation $b > 1, s <1,$ let us introduce the class of L\'evy perpetuities
$$\ELP\; =\; \lacc I(Z) = \int_0^\infty\!\! e^{-Z_s} \,ds, \;\; Z\, =\,\{Z_t, \,t\ge 0\} \; \mbox{is a L\'evy process  with}\; Z_t\to +\infty \;\mbox{a.s.}\racc.$$
Recall that $I(Z)$ is an a.s. convergent integral if and only if $Z$ drifts towards $+\infty$ - see Theorem 1 in \cite{BYS}. Observe also that $I(Z)$ has bounded support when $Z$ is a subordinator with positive drift, and hence may not be in $\II.$ However, more can be said when $Z$ is spectrally negative ($Z$ is an SNLP for short). Introduce the subclass 
$$\ELP_-\; =\; \lacc I(Z) = \int_0^\infty\!\! e^{-Z_s}\, ds, \;\; Z\, =\,\{Z_t, \,t\ge 0\} \; \mbox{is an SNLP with positive mean}\racc.$$
Recall - see the introduction to Chapter 7 in \cite{B} - that an SNLP is characterized by its moment generating function, whose logarithm reads
$$\log \EE[e^{\lbd Z_1}]\; =\; a\lbd \, +\, b \lbd^2\, +\, \int_{-\infty}^0 (e^{\lbd x} - 1 - \lbd x)\, \nu(dx), \quad \lbd\, \ge 0,$$
with $a \in \R, b\in\R^+$ and $\nu$ a non-negative measure on $(-\infty,0)$ whose integral along $1\wedge x^2$ is finite, and that an SNLP drifts towards $+\infty$ if and only if it has a positive mean. The fact that $\ELP_-\subset\SS$ follows from the Markov property at the a.s. finite stopping time $T_x = \inf\{t > 0, Z_t =x\}$ with $x > 0,$ which implies the decomposition
$$I(Z) \; \elaw\; e^{-x} I(Z)\; +\; \int_0^{T_x}e^{-{\tilde Z}_s} ds,$$
where ${\tilde Z}$ is an independent copy of $Z.$ This observation, which is folklore and may be traced back to Vervaat \cite{V} in a more general framework, was used in \cite{NOUS} to prove the self-decomposability of $\ga_a^{-s}, s < 1.$

\begin{THEC} One has
$$\B_{a,b}^{-s}\,\in\, \ELP_-\;\Longleftrightarrow\; b\wedge s \, \le \, 1\,\le\, 2a + b + s +bs.$$
\end{THEC}

\bigskip

\begin{minipage}{0.4\linewidth}
\centering
\psset{unit=1.2cm,algebraic=true,linewidth=0.6pt,
,arrowsize=4pt 2,arrowinset=0.2}

\begin{pspicture*}(-1,-1)(4,4)
\rput[t](3.6,-0.28){$s$}
\rput[r](-0.2,3.6){$b$}

\pspolygon[hatchcolor=couleur,linecolor=white,
fillstyle=hlines,hatchangle=45,hatchsep=0.1](0,0)(0,5)(1,5)(1,1)(5,1)(5,0)
\psline[linecolor=couleur](1,1)(1,10)
\psline[linecolor=couleur](1,1)(10,1)

\pscustom[linecolor=white,fillcolor=white,fillstyle=solid]{\psplot{0}{0.8}{(0.8-x)/(1+x)}\lineto(1,0)\lineto(0,0)\closepath}
\psplot[linecolor=couleur]{0}{0.8}{(0.8-x)/(1+x)}

\psaxes[Dx=1,Dy=1,ticksize=-2pt 2pt]{->}(0,0)(-0.9,-0.9)(4,4)

\end{pspicture*}\\
The case $a<1/2.$
\end{minipage}
\begin{minipage}{0.05\linewidth}
$ $
\end{minipage}
\begin{minipage}{0.4\linewidth}
\centering
\psset{unit=1.2cm,algebraic=true,linewidth=0.6pt,
,arrowsize=4pt 2,arrowinset=0.2}

\begin{pspicture*}(-1,-1)(4,4)
\rput[t](3.6,-0.28){$s$}
\rput[r](-0.2,3.6){$b$}

\pspolygon[hatchcolor=couleur,linecolor=white,
fillstyle=hlines,hatchangle=45,hatchsep=0.1](0,0)(0,5)(1,5)(1,1)(5,1)(5,0)
\psline[linecolor=couleur](1,1)(1,10)
\psline[linecolor=couleur](1,1)(10,1)


\psaxes[Dx=1,Dy=1,ticksize=-2pt 2pt]{->}(0,0)(-0.9,-0.9)(4,4)
\end{pspicture*}\\
The case {$a\ge1/2.$}
\end{minipage}
\vspace{0.6cm}

Combining the three above theorem and a simple asymptotic analysis which will be performed at the end of the next section, we can deduce the main result of the present paper.

\begin{COR} One has $\Babs\in\II$ for all $a,b,s > 0.$ Moreover, one has $\Babs\in\SS\,\Leftrightarrow\, 2a +b +s +bs \ge 1.$ 
\end{COR}

Notice that by (\ref{logB}), this corollary allows to recovers the characterization of $-\log(\Bab)\in\SS$ which is known to be $2a + b \ge 1$ - see again Example VI.12.21 in \cite{SVH}. Towards the end of this paper we will observe that there are situations where $\Babs\in\SS\cap\GG^c.$ It is also worth mentioning that the set of parameters where $\Babs\in\ELP_{-}\cap\HH,$ characterized by the cases (2) and (3) of Theorem 2, is thicker than the set of parameters for $\ga_a^{-s}\in\ELP_{-}\cap\HH$ which is the only line $\{s=1\}$ - see Remark 1 (c) and Section~3.2 in \cite{NOUS}.

Our previous paper \cite{NOUS} had left unanswered the question whether $\ga_a^{-s}\in\GG$ for  $s<1.$ This problem was motivated by the fact that $\ga_a^{-s}\in\HH$ for $s\ge 1$ but not for $s < 1$ since otherwise $\ga_a^{s}$, which is not infinitely divisible, would be also in $\HH$ - see \cite{NOUS} for details and references. This question is also mentioned as an open problem in Section~9 (vi) of \cite{BJTP}. Our last result provides a positive answer.
 
\begin{THED} One has $\ga_a^{-s}\in\GG$ for all $a, s > 0.$
\end{THED}

The proof of this theorem is actually a rather simple consequence of the main result of \cite{BJTP}, which states the important property that the class $\GG$ is stable by independent multiplication. To conclude this introduction, let us also give a few words about the proofs of the three other theorems. For Theorem 1, we use Steutel's characterization of mixtures of exponentials and a logarithmic transformation. The proof of Theorem 2, which is not as immediate as for the HCM characterization of $\ga_a^{-s},$ relies upon Stieltjes transforms and the maximum principle for harmonic functions. For Theorem 3 we  shall use Bertoin-Yor's characterization of $\ELP_-$, as well as several properties of the classical hypergeometric series, old and recent. All these proofs are given in the next section. We conclude the paper with several remarks.

\section{Proofs}

\subsection{Proof of Theorem 1} The density of $\beta_{a,b}^{-s}-1$ reads 
$$f_{a,b,s} (x)\; =\; \frac{\Ga(a+b)}{s\Ga(a)\Ga(b)}\,(x+1)^{\frac{1-a-b}{s} -1}((x+1)^{\frac{1}{s}} -1)^{b-1}\Un_{(0,+\infty)}(x)$$
and is not log-convex if $b > 1$ because 
$$(\log f_{a,b,s})' (x)\; =\; \frac{1}{s(x+1)}\lp\,\frac{b-1}{(x+1)^{\frac{1}{s}} -1}\, -\, a\, -\, s\rp\; \sim\; \frac{b-1}{x}\quad \mbox{as $x\to 0.$}$$
By the Schwarz inequality - see the proof of Theorem 51.6 in \cite{S}, this shows the only if part. The same computation shows that $f_{a,b,s}$ is log-convex for $b\le 1,$ which is known to be sufficient for infinite divisibility - see Theorem 51.4 in \cite{S}. However, it does not seem easy to show directly the reinforcement that $f_{a,b,s}$ is actually CM whenever $b\le 1.$ 

To do so, we first observe that by Corollary 3 in \cite{JS}, it is enough to show that $-\log(\beta_{a,b})\in\MM.$ This latter property follows from Steutel's theorem as had already been noticed in Example VI.12.21 in \cite{SVH}, but we will sketch the argument for the sake of completeness. Using Malmst\'en's formula for the Gamma function - see e.g. Formula 1.9(1) p. 21 in \cite{EMOT}, we first compute
\begin{eqnarray*}
\EE[e^{-\lbd(-\log(\beta_{a,b})})]\;=\; \EE[\beta_{a,b}^\la] &= & \frac{\Ga(a+\la)}{\Ga(a)}\times\frac{\Ga(a+b)}{\Ga(a+b+\la)}\\
&=&\ex{-\int_0^\infty (1-e^{-\la x})\,\frac{e^{-ax}-e^{-(a+b)x}}{x(1-e^{-x})}\,dx}.
\end{eqnarray*}
A further computation similar to the proof of Lemma 2 in \cite{JS} shows that
\begin{eqnarray*}
\frac{e^{-ax}-e^{-(a+b)x}}{x(1-e^{-x})} & = & \int_0^\infty\sharp\{n\in\NN, \; a+n \le t < a+b +n\}\, e^{-xt} \,dt, \qquad x >0, 
\end{eqnarray*}
and the counting function inside the integral clearly satisfies the requirement of Theorem 51.12 in \cite{S} (Steutel's theorem) if and only if $b\le 1$. This completes the proof.

\qed

\subsection{Proof of Theorem 2} \label{Section22}Introduce the function
$$g_{a,b,s}(x) \;=\; \frac{s\Gamma(a)\Gamma(b)}{\Gamma(a+b)}\,\lp \frac{x}{s}\rp^{1-b}\,f_{a,b,s} (x).$$
Since the density of  $\beta_{a,b}^{-s}-1$ is $f_{a,b,s} (x)$ and since $g_{a,b,s} (0+) =1,$ it follows from Property (iv) p.68 and Theorem 5.4.1 in \cite{Bd} that $\beta_{a,b}^{-s}-1\, \in\HH$ if and only if $G_{a,b,s} = -\log(g_{a,b,s})$ is a Thorin-Bernstein function. By Theorem 8.2 (ii) in \cite{SSV}, this is equivalent to $G_{a,b,s}'$ being a Stieltjes transform in the sense of Definition 2.1 in \cite{SSV}. We will use the characterization of Stieltjes transforms given by Corollary 7.4 in \cite{SSV}. Compute
$$G_{a,b,s}'(x)\;=\;\frac{a+s}{s(x+1)}\,+\,(b-1)\lp\frac{1}{x}\,-\,\frac{1}{s(x+1)((x+1)^{1/s}-1)}\rp,$$
which is clearly a Stieltjes transform if $b= 1$ or $s=1.$ We henceforth suppose $b\neq 1$ and $s\neq 1.$ If $s < 1/2,$ then $G_{a,b,s}'$ has at least two poles at $e^{\pm2\ii\pi s} -1$ and hence cannot be a Stieltjes transform. On the other hand, if $s \ge 1/2$ a computation shows that $G_{a,b,s}'$ has an analytic continuation on $\C\setminus (-\infty,-1],$ with
$$G_{a,b,s}'(0)\; =\; \frac{s+a}{s}\;+\;\frac{(b-1)(s+1)}{2s}\cdot$$
We henceforth suppose $s \ge 1/2.$ Setting $x = re^{\ii\theta} - 1$ with $r > 0$ and $\theta\in (0,\pi),$ compute
\begin{eqnarray*}
\Im(G_{a,b,s}'(x)) & = & \frac{-(a+s)\sin(\theta)}{sr}\\
& & \quad +\;(1-b)\lp \frac{r\sin(\theta)}{|r e^{\ii\theta}-1|^2} +\frac{\sin(\theta)}{sr|r^{1/s}e^{\ii\theta/s}-1|^2}-\frac{\sin((1+1/s)\theta)}{sr^{1-1/s}|r^{1/s}e^{\ii\theta/s}-1|^2}\rp.
\end{eqnarray*}
Letting $r\to 0$ and $\theta\to\theta_0\in[0,\pi)$ we obtain
$$\Im(G_{a,b,s}'(x))\; =\; \frac{(1-(a+b+s))\sin(\theta)}{sr}\, +\, o(\sin(\theta) r^{-1}),$$
which shows the necessity of the condition $a+b+s \ge 1.$ Letting now $r\to 0$ and $\theta\to\pi$ we obtain
$$\Im(G_{a,b,s}'(x))\; =\; \frac{(1-b)\sin(\pi/s)}{sr^{1-1/s}}\, +\, \frac{(1-(a+b+s))\sin(\theta)}{sr}\, +\, o(\sin(\theta) r^{-1}),$$ 
which shows the necessity of the condition $(b-1) \sin(\pi/s) \ge 0.$ All in all we have shown that $\beta_{a,b}^{-s}-1\in\HH$ only if (1), (2) or (3) is satisfied, and to finish the proof it remains to prove that $G_{a,b,s}'$ is a Stieltjes transform whenever (1) or (3) holds. To do so, we will use the same argument as in the main theorem of \cite{Bo}. The function 
$$\Im (G_{a,b,s}')\; =\; \Im (g_{a,b,s}'g_{a,b,s}^{-1})$$ 
is harmonic on the open upper half-plane as the imaginary part of an analytic function. Since $G_{a,b,s}'(x)\to 0$ as $\vert x\vert\to +\infty,$ by the maximum principle it is enough to show that 
$$\limsup_{\Im(x)\downarrow 0+} \lp \Im (G_{a,b,s}')(x)\rp\; \le \; 0.$$
Set $l = \lim_{\Im(x)\downarrow 0+} \Re(x+1).$ If $l > 0$ the above limit is clearly zero. If $l < 0,$ a computation shows that the limit equals
$$\frac{(1-b) \sin(\pi/s) \vert l\vert^{1/s-1}}{s \vert \vert l\vert^{1/s} e^{\ii\pi/s} -1\vert^2}$$
and is non-positive if (1) or (3) holds. Last, if $l= 0$ the above analysis for the only if part shows that the limsup is also non-positive.
 
\qed

\subsection{Proof of Theorem 3} We will use the following criterion, which is an immediate consequence of the proof of Proposition 2 in \cite{BYT}.

\begin{FCT}[\bf Bertoin-Yor] Let $X$ be a positive random variable. Then $X\in\ELP_{-}$ if and only $1/X$ has finite exponential moments and there exists an {\rm SNLP} with positive mean whose Laplace exponent $\Psi$ is such that
$$\Psi (u) \; =\; \frac{u  \EE[X^{-(u+1)}]}{\EE[X^{-u}]}$$
for all $u > 0.$
\end{FCT}

Notice that since $1/\Babs = \B_{a,b}^s$ has bounded support, the first condition is always fulfilled. Computing
$$\frac{u  \EE[\B_{a,b}^{s(u+1)}]}{\EE[\B_{a,b}^{su}]} \; = \;\frac{u\Ga(a+s+su)\Ga(a+b+su)}{\Ga(a+b+s+su)\Ga(a+su)}$$
and applying the above criterion, we decuce that $\Babs\in\ELP_-$ if and only if
$$\Psi : u\; \mapsto\; u\times\frac{\Ga(a+s+u)}{\Ga(a+b+s+u)}\times\frac{\Ga(a+b+u)}{\Ga(a+u)}$$
is the Laplace exponent of an {\rm SNLP} with positive mean. Observe that the expression of $\Psi(u)$ is symmetric in $b$ and $s$. Using Gauss' summation formula - see e.g. Formula 2.1.3.(14) p.61 in \cite{EMOT}, we obtain the two transformations
$$\Psi(u) \;= \;u\times\,\! _2F_1(b,-s;a+b+u;1)\; =\; u\times\!\, _2F_1(s,-b;a+s+u;1)$$
where
$$_2F_1(\lbd, \mu;\nu;z)\;=\;\sum_{n\ge0}\frac{(\lbd)_n(\mu)_n}{(\nu)_n\, n!}\,z^n$$
is the classical hypergeometric series, which is defined via the Pochammer symbols $(x)_0 = 1$ and $(x)_n=x(x+1)\dots(x+n-1)$ if $n\ge 1.$ These transformations imply
\begin{equation}
\label{TwoTF}
\Psi(u) \;= \;u\,-\, bsu\sum_{n\ge 1}\frac{(1+b)_{n-1}(1-s)_{n-1}}{(a+b+u)_n n!} \;= \;u\,-\, bsu\sum_{n\ge 1}\frac{(1+s)_{n-1}(1-b)_{n-1}}{(a+s+u)_n n!}\cdot
\end{equation}
Using the first equality in (\ref{TwoTF}) and the partial fraction decomposition
$$\frac{(n-1)!}{(a+b+u)_n}\; =\; \sum_{k=0}^{n-1}\binom{n-1}k \frac{(-1)^k}{a+b+k+u}\; =\; \int_0^\infty e^{-(a+b+u)x} (1-e^{-x})^{n-1} dx,$$
we deduce from Fubini's theorem 
\begin{eqnarray*}
\Psi(u) & = & u\,-\, bsu \int_0^\infty e^{-(a+b+u)x}\lp \sum_{n\ge 1}\frac{(1+b)_{n-1}(1-s)_{n-1}}{n! (n-1)!} (1-e^{-x})^{n-1}\rp dx \\
& = & u\,-\, u \int_0^\infty \!\!e^{-ux}\, \rho_{a,b,s} (x)\, dx
\end{eqnarray*}
with the notation
$$\rho_{a,b,s} (x)\; =\; bs e^{-(a+b)x}\, _2F_1(1+b,1-s;2;1-e^{-x}).$$
Similarly, the second equality in (\ref{TwoTF}) entails
$$\Psi(u)\; =\; u\,-\, u \int_0^\infty \!\!e^{-ux}\, \rho_{a,s,b} (x)\, dx,$$
from which we deduce $\rho_{a,b,s} = \rho_{a,s,b}$ by Laplace inversion. The latter identity, in accordance with the initial symmetry in $(b,s)$ of our problem, can also be obtained from one of Kummer's formul\ae\, for the hypergeometric function - see e.g. 2.1.4(23) p.64 in \cite{EMOT}. In the following, we hence may and will suppose that $s \ge b.$ If $s > b,$ it follows again from Gauss' summation formula that
$$_2F_1(1+b,1-s;2;1-e^{-x})\; \to\; \frac{\Ga(s -b)}{\Ga(1+s)\Ga(1-b)}, \qquad x \to +\infty$$
and this formula extends by continuity to the case $s=b\in\NN,$ the right-hand side being replaced by 
$$\frac{(-1)^b(b-1)!}{\Ga(1+s)}\cdot$$
Last, if $s=b\not\in\NN,$ the asymptotic expansion 2.3.1(2) p.74 in \cite{EMOT} yields 
$$_2F_1(1+b,1-b;2;1-e^{-x})\; \sim\; \frac{x}{\Ga(1+b)\Ga(1-b)}, \qquad x \to +\infty.$$  
In all cases, we deduce that $\rho_{a,b,s}\to 0$ as $x\to +\infty$ and since $\rho_{a,b,s}$ is clearly smooth on $[0, +\infty),$ an integration by parts entails finally
$$\Psi(u)\; =\; u\,+\, \int_0^{+\infty} (1-e^{-ux})\, \rho'_{a,s,b} (x)\, dx\; =\; u\,-\, \int_{-\infty}^0 (e^{ux}-1)\, \rho'_{a,s,b} (-x)\, dx, \quad u\ge 0.$$
It is now clear from the above that $\Psi$ is the Laplace exponent of an SNLP if and only if $\rho'_{a,s,b}$ is non-positive on $[0, +\infty).$ We first compute
$$\rho'_{a,s,b}(0)\; =\; -\frac{bs}2(2a+b+bs+s-1),$$
which shows the necessity of the condition $2a+b+bs+s\ge1.$ Supposing $b < 1,$ we have $(s+1)(1-b) > 0$ and a direct consequence of Theorem 1.3(2) in \cite{AVV} is the log-concavity of $\rho_{a,s,b}.$ In particular, this function is non-increasing on $[0, +\infty)$ if $\rho'_{a,s,b}(0) \le 0,$ which shows the sufficiency of $2a+b+bs+s\ge1.$ If $b =1,$ a simple computation yields 
$$\rho_{a,s,1}(x)\; =\; s e^{-(a+s)x},$$ 
a decreasing function. All in all, we have shown that $\Psi$ is the Laplace exponent of an SNLP whenever $b = b \wedge s\le 1\le 2a + b + s +bs,$ and to conclude the characterization of Theorem B it is enough to show that $\rho'_{a,s,b}$ takes positive values when $b = b \wedge s > 1.$

\medskip

Suppose first that $b\ge 2$ is an integer. Formula 10.8(16) p.170 in \cite{EMOT} shows that
$$_2F_1(1+s,1-b;2;1-e^{-x})\; =\; \frac{1}{n+1} \, P^{(\a, \beta)}_n (2e^{-x}-1)$$
where $P^{(\a, \beta)}_n$ is a Jacobi polynomial with parameters $\a = 1, \beta = s-b \ge 0$ and degree $n = b-1 \in\NN^*.$ By Section~10.16 p.202 in \cite{EMOT}, the latter has $n$ simple roots on $(-1,1)$ and since $\rho_{a,b,s}\to 0$ at infinity, this clearly entails that $\rho'_{a,s,b}$ takes positive values. If $b > 1$ is not an integer, a theorem of F. Klein - see Formel (4), (18) and Satz p.587 in \cite{K} shows that $x\mapsto _2F_1(1+b,1-s;2;1-e^{-x})$
vanishes $[b]$ times on $(0, +\infty),$ where $[x]$ denotes the integer part of a real number $x$. As above, this entails that $\rho'_{a,s,b}$ takes positive values. 

\qed

\subsection{Proof of Theorem 4} Setting $X_{a,s}=-s\log(\ga_a)$, we first compute its Laplace transform
\begin{eqnarray*}
\EE\lcr e^{-\la X_{a,s}}\rcr\; =\;\E{\ga_a^{s\la}} & = &\frac{\Ga(a+s\la)}{\Ga(a)}\;
= \; \ex{\int_0^\infty \lp s\la - \frac{e^{(1-a)x} (1-e^{-s\la x})}{1-e^{-x}}\rp\frac{e^{-x}}{x}dx}, \quad\la \ge 0,
\end{eqnarray*}
where the last equality follows at once from the aforementioned Malmst\'en's formula for the Gamma function. After a change of variable, we can rearrange this expression into
$$\EE\lcr e^{-\la X_{a,s}}\rcr\; =\; \ex{c_{a,s}\la-\int_0^\infty (1-e^{-\la x}-\la x)\,\frac{m_{a,s} (x)}{x}\,dx}$$
with 
$$c_{a,s}\; =\; s\int_0^\infty \lp \frac{e^{-x}}{x} - \frac{e^{-ax}}{1-e^{-x}}\rp dx\qquad \mbox{and}\qquad m_{a,s} (x) \; =\; \frac{e^{-ax/s}}{(1-e^{-x/s})}\cdot$$
For all $\e > 0,$ introduce the random variable $X_{a,s,t}$ with Laplace transform
$$\EE\lcr e^{-\la X_{a,s,t}}\rcr\; =\; \ex{c_{a,s}\la-\int_0^\infty (1-e^{-\la x}-\la x)\,\frac{m_{a,s} (x+\e)}{x}\,dx}, \qquad \lbd \ge 0.$$
Let us rewrite
$$\EE\lcr e^{-\la X_{a,s,t}}\rcr\; =\; \ex{c_{a,s,\e}\la-\int_0^\infty (1-e^{-\la x})\,\frac{m_{a,s} (x+\e)}{x}\,dx}$$
with
$$c_{a,s,\e}\;=\; c_{a,s}\;+\;\int_\e^\infty m_{a,s}(x)\,dx.$$
Since $x \mapsto m_{a,s} (x+\e)$ is CM on $(0,+\infty),$ the random variable $X_{a,s,\e} - c_{a,s,\e}$ belongs to $\GG$ and, by Theorem 3 in \cite{BJTP}, so does the random variable $\rho_{a,s,\e}e^{X_{a,s,\e}} - 1$, with $\rho_{a,s,\e}=e^{-c_{a,s,\e}} > 0.$ This clearly entails that $e^{X_{a,s,\e}}\in\GG$ as well. Letting $\e \to 0$ and using the fact that the class $\GG$ is stable under convergence in law we finally obtain  
$$e^{X_{a,s,\e}}\,\claw\, e^{X_{a,s}}\,\elaw\, \ga_a^{-s}\, \in \, \GG.$$

\qed

\subsection{Proof of the Corollary} By Theorems 1-3, it is clearly enough to show that 
$\beta_{a,b}^{-s}-1\not\in\SS$ as soon as $2a + b + s + bs < 1.$ Supposing $\beta_{a,b}^{-s}-1\in\SS,$ there exists a non-increasing function $k_{a,b,s}$ such that $$\Phi(\la)=\EE\lp e^{-\la (\beta_{a,b}^{-s}-1)}\rp=\exp\lp{-\int_0^\infty (1-e^{-\la x})\frac{k_{a,b,s}(x)}xdx}\rp.$$
With the notation of Section~\ref{Section22}, we compute
\[
-\frac{\Phi'(\la)}{\Phi(\la)} = \frac{{\displaystyle \int_0^\infty e^{-\la x} x^b g_{a,b,s}(x)dx}}{{\displaystyle \int_0^\infty e^{-\la x} x^{b-1} g_{a,b,s}(x)dx} }\cdot
\]
Noticing that $g_{a,b,s}(x) = 1+c_{a,b,s}\,x+o(x)$ as $x\to0$ with $c_{a,b,s}=(1-2a-b-s-bs)/2s$, for every $\al>0$ we obtain by dominated convergence
$$\int_0^\infty e^{-\la x}x^\al g_{a,b,s}(x)dx=\frac1{\la^{\al+1}}\lp \Ga(\al+1) + \frac{c_{a,b,s}\,\Ga(\al+2)}\la  + o(\la^{-1})\rp\quad\text{as $\la\to+\infty$}.$$
This shows the asymptotic expansion
\[
-\frac{\Phi'(\la)}{\Phi(\la)}\; =\; \frac b\la\lp 1+\frac{c_{a,b,s}}\la+o(\la^{-1})\rp\quad\text{as $\la\to+\infty,$}
\]
and since 
$$\la \frac{\Phi'(\la)}{\Phi(\la)}\, +\, b\;=\;\int_0^\infty e^{- x}\lp b- k_{a,b,s}\lp\frac x\la\rp \rp dx,$$
we first deduce $k_{a,b,s}(0+) = b$. In the same way, writing
\[
\la\lp \la \frac{\Phi'(\la)}{\Phi(\la)}+b\rp\; =\; \int_0^\infty xe^{- x}\lp b- k_{a,b,s}\lp\frac x\la\rp \rp \frac\la x\,dx\; \xrightarrow[\la\to+\infty]{}\;b\,c_{a,b,s}.
\] 
shows that
\[
k'_{a,b,s}(0+) \; = \; \frac{b(1-(2a+b+s+bs))}{2s}
\]
and that $k_{a,b,s}$ is a non-increasing function only if $2a+b+s+bs\ge1$.

\qed

\section{Further remarks}

\subsection{On the $\GG-$property for $\beta_{a,b}^{-s}$} Let us first observe that contrary to those of the Gamma distribution, the negative powers of the Beta distribution are not always in $\GG.$ Indeed, if $\beta_{a,b}^{-s}\in\GG$ then so does $\beta_{a,b}^{-s}-1,$ whose density is 
$$f_{a,b,s}(x) \;\sim\; \frac{\Gamma(a+b)}{s\Gamma(a)\Gamma(b)}\,\lp \frac{x}{s}\rp^{b-1}\qquad\mbox{as $x\to 0+$}.$$ 
Hence, by Theorem 4.1.4 in \cite{Bd} - see also the above proof of the Corollary - the Thorin mass of $\beta_{a,b}^{-s}-1$ is then equal to $b$ and, by Theorem 4.1.1 in \cite{Bd}, the function $g_{a,b,s}$ introduced during the proof of Theorem 2 must be CM. However, a computation shows the first order expansion $(\log (g_{a,b,s}))'(x) = C_1 + C_2 x + o(x)$ at zero, with 
$$C_1\; =\; \frac{1}{2}\lp \frac{1-2a}{s} -1\rp \, -\, \frac{b}{2}\lp \frac{1}{s} +1\rp\qquad\mbox{and}\qquad C_2\; =\; 1\, +\,\frac{a}{s} \, +\, \frac{(b-1)}{12}\lp \frac{1}{s} +1\rp\lp \frac{1}{s} +5\rp.$$
The coefficient $C_2$ is positive for $b\vee s \ge 1,$ but for all $a > 0$ there are some $b,s\in (0,1)$ such that $C_2 <0.$ In these cases the function $g_{a,b,s}$ is not log-convex and hence not CM, so that $\beta_{a,b}^{-s}\not\in\GG.$ Observe finally that these $b,s$ can be chosen such that $2a + b + s + bs \ge 1.$ 

It does not seem easy to characterize the $\GG-$property for $\beta_{a,b}^{-s}$ as we did for $\MM, \HH, \ELP_-$ and $\SS.$ We believe that $\beta_{a,b}^{-s}\in\GG$ at least when $b\vee s \ge 1,$ and a first result in this direction is the following proposition. Notice in passing that this result provides an alternative proof of Theorem 4 thanks to (\ref{BG}) taken along integers.\\

\begin{PROP} One has $\beta_{a,b}^{-s}\in\GG$ for all $a,s > 0$ and $b\in\NN.$ 
\end{PROP}

\proof
By Theorem 3 in \cite{BJTP}, it is enough to show that $-\log(\beta_{a,b})\in\GG.$ Recalling
$$\EE[\beta_{a,b}^\la]\;=\;\ex{-\int_0^\infty (1-e^{-\la x})\,\frac{e^{-ax}-e^{-(a+b)x}}{x(1-e^{-x})}\,dx},$$
we see that the $\GG$ property for $-\log(\beta_{a,b})$ is equivalent to the complete monotonicity of the function
$$x\; \mapsto\;\frac{e^{-ax}-e^{-(a+b)x}}{1-e^{-x}},$$
which is easily characterized by $b\in \NN.$

\endproof

\subsection{On the SNLP's associated with $\beta_{a,b}^{-s}$} Under the condition of Theorem 3, it follows from Proposition 2 in \cite{BYT} that $\Babs$ is the perpetuity of a certain compound Poisson process with unit drift and negative jumps. More precisely, one has
\begin{equation}
\label{Blaw}
\Babs\;\elaw\; \int_0^\infty e^{-(t - N^{a,b,s}_t)}\, dt
\end{equation}
where, using the notation of Section~I.1 p. 12 in \cite{B}, the process $\{N^{a,b,s}_t\!, \, t\ge 0\}$ is compound Poisson with L\'evy measure $\nu_{a,b,s} (dx) =-s^{-2}\rho_{a,b,s}'(s^{-1}x) dx$ on $\R^+.$ If $s = 1,$ we have 
$$\nu_{a,b,1} (dx)\; =\; b(a+b)e^{-(a+b)x} dx$$
and this is a slight extension of Example 4 p. 36 in \cite{BYT}, erroneously attributed to Gjessing and Paulsen (we could not locate this special case in \cite{GP}). If $b =1,$ we have 
$$\nu_{a,1,s} (dx)\; =\; (1+a/s)e^{-(1+a/s)x} dx,$$
providing a representation, which we could not find in the literature, of each negative power $\U^{-c}$ of the uniform law $\U$ on $(0,1)$ as the perpetuity of a compound Poisson process with unit drift and exponential L\'evy measure $(1+1/c)e^{(1+1/c)x}$ on $\R^-.$  Notice that $\nu_{a,1,s}$ has a mass always equals to one (and a mean strictly less than one, in accordance with the positive mean of the SNLP whose perpetuity is $\U^{-s/a}$), whereas $\nu_{a,b,1}$ is not a probability when $b\neq 1.$  \\

When $b$ or $s$ is an integer, the hypergeometric series defining $\rho_{a,b,s}$ is a polynomial and the L\'evy measure $\nu_{a,b,s}$ takes a simpler form. For every $n\ge 2,$ the density of $\nu_{a,n,s}$ on $\R^+$ equals
$$\sum_{k=0}^{n-1} \,(1+(a+k)/s)\,c_{k,n,s}\, e^{-(1+(a+k)/s)x}$$
with 
$$c_{k,n,s} \; =\; \prod_{\tiny \begin{array}{c}p=0\\p\ne k\end{array}}^{n-1}(1- s/(p-k)),$$
and it can be checked directly that this density is non-negative if and only if $s\le 1.$ Similarly and with the same notation, the density of $\nu_{a,b,n}$ on $\R^+$ equals
$$\sum_{k=0}^{n-1} \, b n^{-2} (a+b+k)\, c_{k,n,b}\, e^{-(a+b+k)x/n}.$$
These two examples compute the perpetuities of certain compound Poisson processes with unit drift and hyper-exponential L\'evy measures on $\R^-.$ \\

Changing the variable in (\ref{Blaw}), we obtain
$$b^{-s}\Babs\;\elaw\; \int_0^\infty e^{-(b^s t - N^{a,b,s}_{b^s t})}\, dt$$
for all $b\ge 1, s\le 1,$ which is the perpetuity of an SNLP with Laplace exponent
$$\Psi_{a,b,s}(\lbd)\; =\;b^s \lbd\lp 1\,-\,s^{-1}\int_0^\infty \!\!e^{-\lbd x}\, \rho_{a,b,s} (s^{-1} x)\, dx\rp, \qquad \lbd \ge 0.$$
If $s = 1$ and $b\to +\infty$ one has immediately $\Psi_{a,b,s}(\lbd)\to\lbd(a+\lbd),$ in accordance with (\ref{BG}) and Dufresne's result - see Example 3 p.36 in \cite{BYT}. If $s < 1$ and $b\to +\infty,$ we can also recover the SNLP with infinite variation whose perpetuity is distributed as $\ga_a^{-s}.$ Indeed, rewriting
$$\Psi_{a,b,s}(\lbd)\; =\; b^s\lp 1 - \int_0^\infty\!\! \rho_{a,b,s} (x)\, dx\rp\lbd\; +\; s^{-1}\lbd \int_0^\infty \!\!(1-e^{-\lbd x})\, b^s\rho_{a,s,b} (s^{-1}x)\, dx,$$
we see on the one hand from Formula 2.3.2(14) p.77 in \cite{EMOT} and dominated convergence that
\begin{eqnarray*}s^{-1}\lbd \int_0^\infty \!\!(1-e^{-\lbd x})\, b^s\rho_{a,s,b} (s^{-1}x)\, dx  & \to & \frac{\lbd}{\Ga (1-s)} \int_0^\infty \!\!(1-e^{-\lbd x})\, \frac{e^{-(1+a/s)x}}{(1-e^{-x/s})^{1+s}}\, dx
\end{eqnarray*}
as $b\to +\infty.$ The right-hand side transforms, after an integration by parts similar to the one used for Theorem 3, into
$$\int_{-\infty}^0 (e^{\lbd x} - 1 - \lbd x)\, \frac{e^{(1+a/s)x}(s + e^{x/s} + a(1-e^{x/s}))}{s\Ga (1-s)(1-e^{x/s})^{2+s}}\, dx.$$
On the other hand, a change of variable and another integration by parts yields
\begin{eqnarray*}
b^s\lp 1 - \int_0^\infty\!\! \rho_{a,b,s} (x)\, dx\rp & = & b^s\lp 1 - bs \int_0^1 (1-z)^{a+b-1} \, _2F_1(1+b,1-s;2;z)\, dz\rp\\
& = & (a+b-1) b^s\, \int_0^1 (1-z)^{a+b-2} \, _2F_1(b,-s;1;z)\, dz\\
& = & (a+b-1) b^s\, \int_0^1 (1-z)^{a+s-1} \, _2F_1(1+s,1-b;1;z)\, dz\\
& = & \frac{(a+b-1) b^s}{a+s}\, \, _2F_1(1+s,1-b;a+s+1;1)\\
&= & \frac{\Ga(a+s) \Ga(a+b) b^s}{\Ga(a)\Ga(a+b+s)},
\end{eqnarray*}
where the third equality follows from Kummer's transformation 2.9(2) p.105 in \cite{EMOT}, the fourth equality from Bateman's integral formula 2.4(2) p.78 in \cite{EMOT}, and the fifth equality from Gauss' summation formula. Alternatively, this identity can also be seen from the proof of Theorem 3 - see the beginning of Section~\ref{Section33} thereafter. Putting everything together, we see that
$$\Psi_{a,b,s}(\lbd)\; \to\; \lp\frac{\Ga(a+s)}{\Ga(a)} \rp\lbd\; +\;\int_{-\infty}^0 (e^{\lbd x} - 1 - \lbd x)\, \frac{e^{(1+a/s)x}(s + e^{x/s} + a(1-e^{x/s}))}{s\Ga (1-s)(1-e^{x/s})^{2+s}}\, dx$$
as $b \to +\infty,$ and the limit on the right-hand side is precisely the Laplace exponent of the Lemma in \cite{NOUS}. \\

When the conditions of Theorem 3 are not fulfilled, that is when $\rho_{a,b,s}'$ takes positive values, it is easily seen that the function $\lbd \mapsto \Psi(\ii\lbd)$ is no more the L\'evy-Khintchine exponent of a L\'evy process. By Proposition 2 in \cite{BYT}, this shows that $\Babs$ is not distributed as the perpetuity of a L\'evy process having positive jumps and finite exponential moments, and we believe that it is not distributed as the perpetuity of any L\'evy process at all. In a different direction, let us recall that L\'evy processes having positive jumps might have infinitely divisible perpetuities. It follows indeed from Theorem 2.1(j) in \cite{GP} that for every $a > 0, b> 1$ one has
$$\Bab^{-1}\; \elaw\; 1\; +\; \int_0^\infty e^{-(N^{a,b}_t-t)}\, dt$$
where $\{N^{a,b}_t\!, \, t\ge 0\}$ is a compound Poisson process with L\'evy measure $(a+b-1)(b-1)e^{(1-b) x}dx$ on $\R^+.$ From this example of \cite{GP} let us observe, with our previous notation, the mysterious identity 
$$ \int_0^\infty e^{-(t - N^{a,b,1}_t)}\, dt \; \elaw\; 1\; +\; \int_0^\infty e^{-(N^{a,b}_t-t)}\, dt.$$

\subsection{Other remarks} \label{Section33}The proof of Theorem 3 shows that for every $a,b,s >0$ the function

\begin{eqnarray*}
\lbd\;\mapsto\;\frac{\Ga(a+b+\lbd)\Ga(a+s+\lbd)}{\Ga(a+\lbd)\Ga(a+b+s+\lbd)} & = & \lp 1 - \int_0^\infty\!\! \rho_{a,b,s} (x)\rp\; +\; \int_0^\infty (1-e^{-\lbd x})\rho_{a,b,s} (x)\, dx \\
& = & \frac{\Ga(a+b)\Ga(a+s)}{\Ga(a)\Ga(a+b+s)}\; +\; \int_0^\infty (1-e^{-\lbd x})\rho_{a,b,s} (x)\, dx
\end{eqnarray*}
is Bernstein (without drift but with an additional murder coefficient) if and only if $\rho_{a,b,s}$ is non-negative, which, by the aforementioned theorem of F. Klein, occurs if and only if $b\wedge s\le 1.$ It follows then from Theorem 3.6 (ii) in \cite{SSV} that the function
$$\lbd\;\mapsto\;\frac{\Ga(a+\lbd)\Ga(a+b+s+\lbd)}{\Ga(a+b+\lbd)\Ga(a+s+\lbd)}$$
is CM. The fact that the latter function is actually CM for all $a,b,s > 0$ is an easy consequence of Malmst\'en's formula for the Gamma function - see Theorem 6 in \cite{BI}. \\

The related random variables of the type
$$\lp \Bab^{-1} -1\rp^s \; \elaw \; \lp \frac{\ga_b}{\ga_a}\rp^s, \qquad a,b,s >0,$$
with an independent quotient on the right-hand side, appear in the literature as generalized Beta random variables of the second kind, or GB2 random variables. Their density function is
$$\frac{\Ga(a+b)}{s\Ga(a)\Ga(b)}\,(x^{\frac{1}{s}}+1)^{-(a+b)} x^{\frac{b}{s} -1}\Un_{(0,+\infty)}(x)$$
and easily seen to be HCM if and only if $s\ge 1.$  Let us mention in passing that this contrasts with the independent product $(\ga_a\times\ga_b)^s,$ which belongs to $\HH$ if $\vert b-a\vert\le 1/2\le s$ - see Section~2.2 in \cite{Bo}. It is also clear that the GB2 random variables do not have negative exponential moments and hence cannot belong to $\ELP_-.$ It would be interesting to know whether the GB2 random variables belong to $\II$ notwithstanding, when $s < 1.$ It can be shown rather easily by the $\ga_2$-criterion - see e.g. Theorem VI.4.5 in \cite{SVH} - that this is indeed the case for $a \le 1/2 \le s \le b = 1-a.$ We postpone the remaining cases to future research.

\bigskip

\noindent
{\bf Acknowledgement.} Part of this work was written during a stay at Dresden of the second author, who wishes to thank Anita Behme and Ren\'e Schilling for their hospitality.

\end{document}